\documentclass[12pt]{article}%
\usepackage[applemac]{inputenc}
\usepackage{amsmath,amssymb,fullpage}
\usepackage{amsfonts}
\usepackage{amsmath,amssymb}
\usepackage[applemac]{inputenc}
\usepackage{amsmath,amssymb,fullpage}
\usepackage{color}
\usepackage{color, colortbl}
\usepackage{amsmath}
\usepackage{amssymb}
\usepackage{graphicx}
\usepackage{tikz}
\usepackage{geometry}
\usetikzlibrary{arrows}
\usepackage{amsfonts}
\usepackage{amsmath,amssymb}
\usepackage[applemac]{inputenc}
\usepackage{amsmath,amssymb,fullpage}
\usepackage{color}
\usepackage{amsmath}
\usepackage{amssymb}
\usepackage{graphicx}
\usepackage{tikz}
\newtheorem{theorem}{Theorem}[section]
\newtheorem{lemma}[theorem]{Lemma}

\newtheorem{definition}[theorem]{Definition}

\newtheorem{example}[theorem]{Example}

\newtheorem{problem}[theorem]{Problem}

\usetikzlibrary{arrows,shapes,automata,petri}

\newcommand{\dproof}{\noindent {Proof.} \quad}
\newcommand{\fproof}{\hfill $\square$ \bigskip}

\numberwithin{equation}{section}

\def\RB{\mathbb{R}}

\def \ind{1\!\!1}

\definecolor{LightCyan}{rgb}{0.88,1,1}

\def\MM{{\mathbb{M}}}
\def\RR{{\mathbb{ R}}}

\def\EE{{\mathbb{ E}}}
\def\D{\mathcal{D}}

\def\1B{\text{1\!\!I}}

\def\<{\langle}
\def\>{\rangle}

\begin{document}

\title{Optimal stopping of conditional McKean-Vlasov jump diffusions}
\author{Nacira Agram$^{1}$ \& Bernt \O ksendal$^{2}$}
\date{28 July 2022}
\maketitle

\footnotetext[1]{Department of Mathematics, KTH Royal Institute of Technology 100 44, Stockholm, Sweden. \newline
Email: nacira@kth.se. Work supported by the Swedish Research Council grant (2020-04697).}

\footnotetext[2]{%
Department of Mathematics, University of Oslo, Norway. 
Email: oksendal@math.uio.no.}

\begin{abstract}
We study the problem of optimal stopping of conditional McKean-Vlasov  (mean-field) stochastic differential equations with jumps (conditional McKean-Vlasov jump diffusions, for short). We obtain sufficient variational inequalities for a function to be the value function of such a problem and for a stopping time to be optimal.\\ To achieve this, we combine the state equation for the conditional McKean-Vlasov equation with the associated stochastic Fokker-Planck equation for the conditional law of the solution of the state. This gives us a Markovian system which can be handled by using a version of the Dynkin formula.\\
We illustrate our result by solving explicitly two optimal stopping problems for conditional McKean-Vlasov jump diffusions. More specifically, we first find the optimal time to sell in a market with common noise and jumps,
and, next, we find the stopping time to quit a project whose state is modelled by a jump diffusion, when the performance functional involves the conditional mean of the state.

\end{abstract}


\section{Introduction}
Given a $d$-dimensional Brownian motion $B=(B_1, B_2, ..., B_d)$, a $k$-dimensional compensated Poisson random measure $\tilde{N}$ on a filtered probability space $(\Omega,\mathcal{F},P, \mathbb{F}=\{\mathcal{F}\}_{t\geq 0})$ and a random variable $Z \in L^2(P)$ that is independent of $\mathbb{F}$, we consider the state dynamics
\small
\begin{align}
X(t) &=Z+\int_0^t\alpha(s,X(s),\mu_s)dt+\beta(s,X(s),\mu_{s})dB(s)\nonumber\\
&+\int_0^t\int_{
\mathbb{R}^d
}\gamma(s,X(s^-),\mu_{s^{-}},\zeta)\widetilde{N}(ds,d\zeta),\label{1.1}
\end{align}
where we denote by $\mu_t=\mathcal{L}(X(t) | \mathcal{F}_t^{(1)})$ the conditional probability distribution of $X(t)$ given the filtration $\mathcal{F}_t^{(1)}$generated by the the first component $B_1(u);u\leq t$ of the Brownian motion up to time $t$.
Loosely speaking, the equation above models a McKean-Vlasov dynamics which is subject to what is called a "common noise" coming from the Brownian motion $B_1(t)$, which is observed and is influencing the dynamics of the system.

Standard optimal stopping problems can be solved by the dynamic programming principle (DPP, for short) and the Hamilton-Jacobi-Bellman (HJB) equation; we refer e.g. El Karoui \cite{E}, Shiryaev \cite{S} for earlier works in the continuous case and to \O ksendal and Sulem \cite{OS} in the setting of jump diffusions. In the non-standard optimal stopping problem with inside information, we refer to Hu and \O ksendal \cite{HO} where a Malliavin calculus approach is used.

The process $X(t)$ given by \eqref{1.1} is not in itself Markovian, so to be able to use the HJB, we extend the system to the process $Y$ defined by 
\begin{equation*}
Y(t)=(s+t,X(t),\mu_t); \quad t\geq 0;\quad Y(0)=(s,Z,\mu_0)=:y,
\end{equation*}
 for some arbitrary starting time $s\geq 0$, with state dynamics given by $X(t)$ and conditional law of the state given by $\mu_t$. This system is Markovian, in virtue of the following equation for the conditional law $\mu_t$.\\
  It is proved in \cite{AO} that $\mu_t$ satisfies the following stochastic partial integro-differential equation (SPIDE, for short) of Fokker-Planck type:
\begin{align*} 
d\mu _{t} =A_0^{*} \mu_t dt + A_1^{*}\mu_t dB_1(t);   \quad \mu_0=\mathcal{L}(X(0))=\mathcal{L}(Z),
\end{align*}
where $A_0^{*}, A_1^{*}$ are integro-differential operators which will be specified later.\\

Given the process $Y(t)$, a profit rate $f$, a bequest function $g$ and a family $\mathcal{T}$ of $\mathbb{F}$-stopping times, we consider the following \emph{McKean-Vlasov common noise optimal stopping problem}:
 \begin{problem}
 Find the value function $\Phi(y)$ and an optimal stopping time $\tau^{*} \in \mathcal{T}$ such that 
 \begin{align*}
  \Phi(y):&=  \sup_{\tau \in \mathcal{T}} \mathbb{E}^{y}\Big[ \int_0^{\tau} f(Y(t))dt  + g(Y(\tau)) \Big]
  = \mathbb{E}^{y}\Big[ \int_0^{\tau^{*}} f(Y(t))dt  + g(Y(\tau^{*})) \Big],
  \end{align*}
 where $\mathbb{E}^{y}$ denotes the expectation given that $Y(0)=y$.
 \end{problem}
 This problem will be formulated more precisely later.
 
When there is no dependence on the conditional law on the coefficients $b,\sigma,f$ and $g$, the McKean-Vlasov optimal stopping problem with common noise reduces to the optimal stopping problem with jump diffusion, see for instance \O ksendal and Sulem \cite{OS}.

There has been a lot of attention to mean-field game (MFG, for short), following the works by  Lasry and Lions \cite{LL} and Huang, Malham\' e and Caines \cite{M}. Furthermore, some recent  papers consider optimal stopping of MFG, in particular they are looking for MFG equilibria. We refer  to Bertucci \cite{B}, Bouveret, Dumitrescu and Tankov \cite{BDT}, Dumitrescu, Leutscher and Tankov \cite{DLT}, Carmona, Delarue and Lacker \cite{CDL}, Nutz \cite{N}. More precisely, given $\mu_t=\mathcal{L}(X(t))$ and an initial real value $x$, they consider the state dynamics
\small 
\begin{align*}
X^{\mu}(t) =x+\int_0^t\alpha(s,X^{\mu}(s),\mu_s)dt+\beta(s,X^{\mu}(s),\mu_{s})dB(s),
\end{align*}
and the following optimal stopping problem:
$$\sup_{\tau} \mathbb{E}\Big[ \int_0^{\tau} f(s,X^{\mu}(s),\mu_s)ds  + g(\tau,X^{\mu}(\tau),\mu(\tau)) \Big].$$
In the above mentioned papers \cite{B, BDT, DLT, CDL, N} the authors assume they have obtained an optimal stopping time $\tau^*(\mu)$ and are looking for a mean-field equilibrium.They find a MFG equilibrium by a fixed point argument. In our paper, we do not consider a MFG setup. On the other hand, we are including common noise and jumps; which is not the case for them.

In the recent paper by Talbi, Touzi and Zhang \cite{TTZ} the mean-field stopping problem for jump diffusions studied in a weak formulation
i.e., in terms of the joint marginal law of the stopped underlying process and the survival process. Moreover, in the criterion function, the terminal cost is deterministic. Our paper also differs from that paper in other ways:\\
(i) First, we are proving sufficient variational inequalities , while in \cite{TTZ} a necessary version is proved (assuming that the value function is smooth).\\ 
(ii) Second, we are considering strong solutions.\\
(iii) Third, in \cite{TTZ} they consider the law of the stopped process, while we consider the conditional law of the unstopped process up to the stopping time and then the conditional law of the state at the stopping time. There is a subtle difference between the two models. Which model to choose depends on what applications one has in mind. 

We highlight that, until the present work, none of the above mentioned works find an optimal stopping time explicitly for mean-field jump diffusions with common noise as we do here.

The organisation of our paper is as follows:\\
Section 2 is devoted to recall the stochastic Fokker-Planck equation for the conditional McKean-Vlasov jump diffusion. In Section 3 we study the optimal stopping problem for McKean-Vlasov jump diffusion with common noise and we prove sufficient variational inequalities for optimality. In Section 4, we apply the result of Section 3 to solve explicitly two optimal stopping problems:  \\
(i) Finding the optimal time to sell in a market with common noise and jumps.\\
(ii) Finding the optimal time to quit a project whose state is modelled by a jump diffusion, when the performance functional involves the conditional mean of the state.


\section {Stochastic Fokker-Planck equation for the conditional McKean-Vlasov jump diffusion}
We start this section by first introducing some notations and definitions that will be used throughout this paper.\\
 A Radon measure on $\mathbb{R}^d$ is a Borel measure which is finite on compact sets, outer regular on all Borel sets and inner regular on all open sets. In particular, all Borel probability measures on $\mathbb{R}^d$ are Radon measures.\\
In the following, we let $\mathbb{M}_0$ be the set of deterministic Radon measures and we let $C_0(\mathbb{R}^d)$ be the uniform closure of the space $C_c(\mathbb{R}^d)$ of continuous functions with compact support. If we equip $\mathbb{M}_0$ with the total variation norm $||\mu||:=|\mu|(\mathbb{R}^d)$, then $\mathbb{M}_0$ becomes a Banach space, and it is the dual of $C_0(\mathbb{R}^d)$. See  Chapter 7 in Folland \cite{F} for more information.

Let $X(t)=X_t \in \RR^d$ be a mean-field stochastic differential equation with jumps, from now on called a \emph{McKean-Vlasov jump diffusion}, of the form (using matrix notation),
\begin{align}\label{MK}
dX(t)& =\alpha(t,X(t),\mu_t)dt+\beta(t,X(t),\mu_{t})dB(t)+\int_{
\mathbb{R}^d
}\gamma(t,X(t^-),\mu_{t^{-}},\zeta)\widetilde{N}(dt,d\zeta), \nonumber\\ 
X(0) &=Z,
\end{align}
where $B\in \RR^m = \RR^{m\times 1}, \widetilde{N} \in \RR^k=\RR^{k\times 1}$ are, respectively, an $m$-dimensional  Brownian motion and a $k$-dimensional compensated Poisson random measure. Here $Z$ is a random variable which is independent of the $\sigma$-algebra generated by $(B,\widetilde{N})$ and such that
$$
\EE[|Z|^2]<\infty.
$$
Define the $\sigma$-algebra $\mathbb{F}=\{\mathcal{F}\}_{t\geq0}$ to be the filtration generated by $Z$ and $(B,\widetilde{N})$.\\
We assume that the coefficients: $\alpha(t,x,\mu):[0,T]\times \mathbb{R}^d \times \mathbb{M}\rightarrow \mathbb{R}^d, \beta(t,x,\mu):[0,T]\times \mathbb{R}^d \times \mathbb{M}\rightarrow \mathbb{R}^{d \times m}$ and $\gamma(t,x,\mu,\zeta):[0,T]\times \mathbb{R}^d\times \mathbb{M}\times \mathbb{R}^d \rightarrow \mathbb{R}^{d \times k}$ are $\mathbb{F}$-predictable for all $x,\mu, \zeta$, 
and that $\alpha, \beta, \gamma$ are continuous with respect to $t$ and $x$ for all $\mu,\zeta$. 

One can easily check that under some assumptions such as the Lipschitz and the linear growth conditions, there exists a unique solution of equation  \eqref{MK}. 

For convenience, we assume that for all $\ell; 1 \leq \ell \leq k,$ the L\'evy measure of $N_{\ell}$, denoted by $\nu_{\ell}$, satisfies the condition
$\int_{\RR} \zeta^2 \nu_{\ell}(d\zeta) < \infty$, which means that $N_{\ell}$ does not have many big jumps (but $N_{\ell}$ may still have infinite total variation near $0$). This assumption allows us to use the $L^2$ version of the It\^o formula for jump diffusion given in Theorem 1.16 in \cite{OS}.\\
We denote by $\mu_t=\mu_t(\omega)=\mu_t(dx,\omega) $ the \emph{conditional law} of $X(t)$ given the filtration $\mathcal{F}_t^{(1)}$ generated by the first component $B_1$ of the $m$-dimensional Brownian motion $B$. More precisely, we consider the following model:
\begin{definition}
We assume that $m \geq 2$ and we fix one of the Brownian motions, say $B_1=B_1(t,\omega)$, with filtration $\{\mathcal{F}_t^{(1)}\}_{t\geq 0}$. We define $\mu_t=\mu_t(\omega, dx)$ to be regular conditional distribution of $X(t)$ given $\mathcal{F}_t^{(1)}$. This means that $\mu_t(\omega,dx) $ is a Borel probability measure on $\RR^d$ for all $t \in [0,T],\omega \in \Omega$ and
\begin{equation} \label{cond}
\int_{\mathbb{R}^n} g(x) \mu_t(dx,\omega)= \EE[g(X(t)) | \mathcal{F}_t^{(1)}](\omega)
\end{equation}
 for all functions $g$ such that $\EE[ |g(X(t)) |] < \infty$.
\end{definition}

 From now on we let $\mathbb{M}$ denote all random measures $\lambda(dx,\omega)$ which are Radon measures with respect to $x$ for each $\omega$.  

Recall that the  stochastic Fokker-Planck differential equation for the conditional distribution of the McKean-Vlasov equations has been studied by Coghi and Gess \cite{CG}. It was subsequently extended to jump diffusions by Agram and \O ksendal \cite{AO}, who obtained a corresponding stochastic Fokker-Planck integr0-differential equation for the conditional law.

\begin{theorem}{(Conditional stochastic Fokker-Planck equation \cite{AO})}\\
Let $X(t)$ be as in \eqref{MK} and let $\mu_t=\mu_t(dx,\omega)$ be the regular conditional distribution of $X(t)$ given $\mathcal{F}_t^{(1)}$. Then $\mu_t$ satisfies the following SPIDE (in the sense of distributions):
 \begin{align} 
d\mu _{t} =A_0^{*} \mu_t dt + A_1^{*}\mu_t dB_1(t);   \quad \mu_0=\mathcal{L}(X(0)),\label{A}
\end{align}
where 
$A_0^{*}, A_1^{*}$ are the integro-differential and the differential operator, respectively,  which are given respectively by:
\begin{align}
A_0^{*}\mu&= -\sum_{j=1}^d D_j [\alpha_j \mu] +\frac{1}{2}\sum_{n,j=1}^d D_{n,j}[(\beta \beta^{T})_{n,j} \mu] \nonumber\\
&+\sum_{\ell=1}^k  \int_{\mathbb{R}}\Big\{\mu^{(\gamma^{(\ell)})}-\mu+\sum_{j=1}^d D_j[\gamma_j^{(\ell)}(s,\cdot,\zeta)\mu]\Big\} \nu_{\ell} \left( d\zeta \right), \label{A0*}
\end{align}
and
\begin{align}
A_1^{*}\mu= - \sum_{j=1}^m D_j[\beta_{1,j} \mu], \label{A1*}
\end{align}
were $\beta^{T}$ denotes the transposed of the $d \times m$ - matrix $\beta=\big[\beta_{k,j}\big]_{1\leq k \leq d,1 \leq j \leq m}$  and $\gamma^{(\ell)}$ is column number $\ell$ of the matrix $\gamma$.

\end{theorem}
In the above $D_j, D_{n,j}$
denote $\frac{\partial }{\partial x_j}$  and  $\frac{\partial^2}{\partial x_n \partial x_j}$ respectively, in the sense of distributions and $\mu^{(\gamma^{(\ell)})}$ denotes the $\gamma^{(\ell)}$-shift of $\mu$.


\section {Sufficient variational inequalities for optimal stopping of conditional McKean-Vlasov jump diffusions}
Consider the following $1 \times d \times 1$-dimensional $[0,\infty) \times \mathbb{R}^d\times  \mathbb{M}$ - valued process $Y(t):=(s+t, X(t),\mu_t)$:
\begin{align}\label{Y2}
dY(t)&=\left[\begin{array}{clcr}
dt \\dX(t)\\d\mu_t
\end{array} \right] =\left[ \begin{array}{c}
1\\ \alpha(Y(t)) \\A_0^{*}\mu_t
\end{array} \right] dt
+\left[\begin{array}{rc}
0_{1\times m} \\ \beta(Y(t)) \\ A_1^{*}\mu_t,0,0 ...,0
\end{array} \right]dB(t)\nonumber\\
&+ \int_{\mathbb{R}^d} \left[ \begin{array}{rc}
0_{1\times k}\\ \gamma(Y(t^-),\zeta) \\0_{1\times k}
\end{array}\right] \widetilde{N}(dt,d\zeta),\quad s\leq t \leq T,
\end{align}
where $X(t)$ and $\mu_t$ satisfy the equations \eqref{MK} and \eqref{A}, respectively. 
Moreover, we have used the shorthand notation
\begin{align*}
\alpha(Y(t))&=\alpha(s+t,X(t),\mu(t))\\
\beta(Y(t))&= \beta(s+t,X(t),\mu(t))\\
\gamma(Y(t^-),\zeta)&=\gamma(s+t,X(t^-),\mu(t^-),\zeta).
\end{align*}
The process $Y(t)$ starts at $y=(s, Z,\mu)$. We shall denote by $\mu$ the initial probability distribution $\mathcal{L}(X(0))$ or the generic value of the conditional law $\mu_t :=\mathcal{L}(X(t) | \mathcal{F}_t^{(1)})$, when there is no ambiguity.

We proceed to formulate precisely our optimal stopping problem for conditional McKean-Vlasov jump diffusions:

Let $S$ be a given open set (the solvency region)  in $[0,\infty) \times \RR^d \times \MM$ and define
\begin{equation}
\tau_S=\tau_S(y,\omega)=\inf\{t>0;Y(t,\omega)\notin S\} \text{ (the bankruptcy time)},
\end{equation}
and define $\mathcal{T}$ to be the set of $\mathbb{F}$-stopping times $\tau \leq \tau_S$.\\
Let  $f: [0,\infty) \times \mathbb{R}^d \times \mathbb{M} \mapsto \mathbb{R}$ (the profit rate) and $g: [0,\infty) \times \mathbb{R}^d \times \mathbb{M} \mapsto \mathbb{R}$ (the bequest function) be given functions, such that
\begin{align*}
\mathbb{E}^{y}\Big[ \int_0^{\tau} |f(Y(t))|dt  + |g(Y(\tau))| \Big] < \infty \text{ for all } \tau\leq \tau_S.
\end{align*}
We introduce the performance functional:
\begin{align*}
&J^{\tau}(y)=\mathbb{E}^{y}\Big[ \int_0^{\tau} f(Y(t))dt  + g(Y(\tau)) \Big], \quad \tau \in \mathcal{T}.
\end{align*}
Assume that $f$ and $g$ satisfies the following conditions:
\begin{itemize}
\item
$\EE^y\Big[{\displaystyle\int\limits_0^{\tau_S}} f^-(Y(t))dt\Big]<\infty$;\qquad for
all $y\in S$,

\item
the family $\{g^-(Y(\tau));\tau$ stopping time, $\tau\leq \tau_S\}$
is uniformly integrable 
 for all $y\in S$.
 \end{itemize}
 Here $f^{-}(y)=\max\{-f(y),0\}$ is the negative part of $f$, and similarly with $g$.\\

We consider the following problem: \\
\begin{problem}{(Optimal stopping problem)} Find $\Phi(y)$ and 
$\tau^\ast\in{\cal T},$ 
such that
\begin{equation}
\Phi(y)=\sup_{\tau\in{\cal T}}J^\tau(y)=J^{\tau^\ast}(y)\;,\qquad \mbox{for $\,\tau\in{\cal T}, y \in S$}.
\end{equation}
\end{problem}
Note that since $J^0(y)=g(y)$ we have
\begin{equation}
\Phi(y)\geq g(y)\qquad\mbox{for all $\,y\in S$}\;.
\end{equation}
We will formulate sufficient variational inequalities for this problem:\\
In the following, if $\varphi=\varphi(s,x,\mu) \in C^{1,2,2}([0,T] \times \mathbb{R}^d \times \mathbb{M})$, then  $D_{\mu} \varphi=\nabla_{\mu} \varphi \in L(\mathbb{M},\RR)$ (the set of bounded linear functionals on $\mathbb{M}$) denotes the 
Fr\' echet derivative (gradient) of $\varphi$ with respect to $\mu \in \mathbb{M}$. Similarly $D_{\mu}^2 \varphi$ denotes the double derivative of $\varphi$ with respect to $\mu$. It is an element of $L(\mathbb{M},L(\mathbb{M},\RR))$. By the Riesz representation theorem it may be regarded as an element of $L(\mathbb{M} \times \mathbb{M},\RR)$ (the bounded linear functionals on $\mathbb{M} \times \mathbb{M}$).\\
The infinitesimal generator $\mathcal{G}$ of the Markov jump diffusion process $Y(t)$ is defined by  
\begin{align*}
&\mathcal{G} \varphi= \frac{\partial \varphi}{\partial s} +\sum_{j=1}^d \alpha_j \frac{\partial \varphi}{\partial x_j} + \langle \nabla_{\mu} \varphi, A_0^{*} \mu \rangle + \tfrac{1}{2}\sum_{j,n=1}^{d}  (\beta \beta^{T})_{j,n}\frac{\partial ^2 \varphi}{\partial x_j \partial x_n} \nonumber\\
& + \tfrac{1}{2}\sum_{j=1}^d \beta_{j,1}\frac{\partial}{\partial x_j}\langle\nabla_{\mu} \varphi,A_1^{*}\mu\rangle +\tfrac{1}{2} 
\langle A_1^{*}\mu, \langle D_{\mu}^2 \varphi,A_1^{*}\mu\rangle \rangle \nonumber\\
&+\sum_{\ell =1}^k \int_{\mathbb{R}} \{ \varphi(s, x+\gamma^{(\ell)}, \mu)) - \varphi(s,x,\mu) 
-\sum_{j=1}^d\gamma_j^{(\ell)}  \tfrac{\partial}{\partial x_j} \varphi(s,x,\mu) \}\nu_{\ell}(d\zeta) ;\nonumber\\
& \varphi  \in C^{1,2,2}([0,T] \times \mathbb{R}^d \times \mathbb{M}).
\end{align*}
Here  $\gamma=\gamma(s,x,\mu,\zeta)$ and $\gamma^{(\ell)}$ is column number $\ell$ of the $d \times k$- matrix $\gamma$.\\

We can now formulate the associated variational inequalities for this problem:
\begin{theorem}\textbf{(Variational inequalities for optimal stopping of conditional McKean-Vlasov jump diffusions)} 
\begin{itemize}
    \item [(A)] Suppose we can find a function $\phi\colon \overline{S}\to\RR$
such that
\begin{description}
\item[{\rm (i)}]
$\phi\in C^1(S)\cap C(\overline{S}\,)$
\item[{\rm (ii)}]
$\phi\geq g$ on $S$ and $\lim\limits_{t\to\tau_S^-} \phi(Y(t))=g(Y({\tau_S}))
{\ind{}}_{\{\tau_S<\infty\}}$ a.s.
\end{description}

Define
$$
\D=\{x\in S;\phi(x)>g(x)\}\;.
$$
Suppose $Y(t)$ spends 0 time on $\partial \D$ a.s., i.e.
\begin{description}

\item[{\rm (iii)}]
$\EE^y\Big[{\displaystyle\int\limits_0^{\tau_S}
{\ind{}}_{\partial \D}(Y(t))dt\Big]=0}$ 
for all $y\in S$
\smallskip

and suppose that
\smallskip

\item[{\rm (iv)}]
$\partial \D$ is a Lipschitz surface, i.e. $\partial \D$ 
is locally the graph of a function $h\colon \RR^{k-1}\to\RR$ such that
there exists $K<\infty$ with
$$
|h(x)-h(y)|\leq K|x-y|,\qquad \mbox{for all $\,x,y$}\;.
$$
\end{description}
Moreover, suppose the following:
\begin{description}

\item[{\rm (v)}]
$\phi\in C^2(S\setminus\partial \D)$ and the second order derivatives
of $\phi$ are locally bounded near $\partial \D$

\item[{\rm (vi)}]
$\mathcal{G} \phi+f\leq0$ on $S\setminus \D$.

Then
$$
\phi(y)\geq \Phi(y),\qquad \hbox{for all $y\in S$}.
$$
\end{description}
\end{itemize}
\begin{itemize}
    \item [(B)] Suppose, in addition to the above, that
\begin{description}
\item[{\rm (vii)}]
$\mathcal{G} \phi+f=0$ on $\D$

\item[{\rm (viii)}]
$\tau_{\D}\colon =\inf\{t>0;Y(t)\notin \D\}<\infty$ a.s. for all  
$y\in S$

and

\item[{\rm (ix)}]
the family $\{\phi(Y(\tau));\tau\leq\tau_{\D},\tau\in{\cal T}\}$ is 
uniformly integrable, for all $y\in S$. Then
\end{description}
\end{itemize}
\vspace*{-\belowdisplayskip}
\begin{equation}
\phi(y)=\Phi(y)=\sup_{\tau\in{\cal T}}\EE^y\bigg[
\int\limits_0^\tau f(Y(t))dt+g(Y(\tau))\bigg]\;;
\qquad y\in S
\end{equation}
\vspace*{-3ex}
and
\begin{equation}
\tau^\ast=\tau_{\D}
\end{equation}

is an optimal stopping time for this problem.\\

\end{theorem}

\dproof
By (i), (iv) and (v) we can find a sequence of functions 
$\phi_j\in C^2(S)\cap C(\overline{S}\,)$,\; \\$j=1,2,\ldots$, such that
\begin{itemize}
    \item [(a)] $\phi_j\to \phi$ uniformly on compact subsets of $\overline{S}$, as
$j\to\infty$,
    \item [(b)] $\mathcal{G}\phi_j\to \mathcal{G}\phi$ uniformly on compact subsets of $S \setminus \partial \D$, as $j\to\infty$,
    \item [(c)] $\{\mathcal{G}\phi_j\}_{j=1}^\infty$ is locally bounded on $S$.
\end{itemize}
See Theorem 3.1 in \cite{OS}.The proof there applies with minor changes.\\
Let $\big\{S_R\big\}_{R=1}^\infty$ be a sequence of bounded open sets
such that $S=\bigcup\limits_{R=1}^\infty S_R$. Put 
$T_R=\min(R,\inf$ $\{t>0;Y(t)\not\in S_R\})$ and let $\tau\leq \tau_S$ be
a  stopping time. Let $y\in S$. Then by Dynkin's formula
\vspace*{-\belowdisplayskip}
\begin{equation}
\EE^y[\phi_j(Y({\tau\wedge T_R}))]=\phi_j(y)
    +\EE^y\bigg[\int\limits_0^{\tau\wedge T_R}\mathcal{G}  \phi_j(Y(t))dt\bigg].
\end{equation}
Hence by (a), (b), (c) and (iii) and bounded a.e. convergence
\begin{eqnarray}
\phi(y)&=&\lim_{j\to\infty} 
    \EE^y\bigg[\int\limits_0^{\tau\wedge T_R}-\mathcal{G}  \phi_j
     (Y(t))dt+\phi_j(Y({\tau\wedge T_R}))\bigg] \nonumber \\
&=&
   \EE^y\bigg[\int\limits_0^{\tau\wedge T_R}-\mathcal{G}  \phi(Y(t))dt+
    \phi(Y({\tau\wedge T_R}))\bigg]\;.
\end{eqnarray}
Therefore, by (ii), (iii) and (vi),
$$
\phi(y)\geq 
\EE^y\bigg[\int\limits_0^{\tau\wedge T_R}f(Y(t))dt+
g(Y({\tau\wedge T_R}))\bigg]\;.
$$
Hence by the Fatou lemma, we get
$$
\phi(y)\geq\lim_{\overline{R\to\infty}} 
\EE^y\bigg[\int\limits_0^{\tau\wedge T_R}
f(Y(t))dt+g(Y({\tau\wedge T_R}))\bigg]\geq
\EE^y\bigg[\int\limits_0^\tau f(Y(t))dt+g(Y(\tau))\bigg]\;.
$$
Since $\tau\leq \tau_S$ was arbitrary, we conclude that
\begin{equation}\label{4.9}
\phi(y)\geq\Phi(y)\qquad\mbox{for all $\,y\in S$}\;,
\end{equation}
which proves (A). \\

We proceed to prove (B):\\
If $y\notin \D$ then $\phi(y)=g(y)\leq\Phi(y)$ so by \eqref{4.9}
we have
\begin{equation}\label{4.10}
\phi(y)=\Phi(y)\quad\mbox{and}\quad
\widehat{\tau}=\widehat{\tau}(y,\omega)\colon =0\qquad\mbox{is
optimal for $\,y\notin \D$}\;.
\end{equation}
Next, suppose $y\in \D$. 
Let $\{\D_k\}_{k=1}^\infty$ be an increasing
sequence of open sets $\D_k$ such that $\overline{\D}_k\!\subset\! \D$,
$\,\overline{\D}_k$ is compact and  
$\D\!=\!\bigcup\limits_{k=1}^\infty\! \D_k$.\\ 
Put
${\tau_k\!=\!\inf\{t\!>\!0; Y(t
)\!\not\in\! \D_k\}}$,  
$k=1,2,\ldots$
By Dynkin's formula we have for $y\in \D_k$,
\begin{align*}
\phi(y) &= \lim_{j\to\infty}\phi_j(y)=
    \lim_{j\to\infty}\EE^y\bigg[\int\limits_0^{\tau_k\wedge T_R}-
    \mathcal{G} \phi_j(Y(t))dt+\phi_j(Y({\tau_k\wedge T_R}))\bigg] \\
&= \EE^y\bigg[\!\int\limits_0^{\tau_k\wedge T_R}\!- \mathcal{G} \phi(Y(t))dt+
     \phi(Y({\tau_k\wedge T_R}))\bigg]\! \\
&=\EE^y\bigg[\!\int\limits_0^{\tau_k\wedge T_R}\! 
    f(Y(t))dt
   +\phi(Y({\tau_k\wedge T_R}))\bigg].
\end{align*}
So by uniform integrability and (ii), (vii), (viii), we get
\begin{eqnarray}\label{4.11}
\phi(y) &=& \lim_{R,k\to\infty}\EE^y\bigg[
     \int\limits_0^{\tau_k\wedge T_R}
    f(Y(t))dt+\phi(Y({\tau_k\wedge T_R}))\bigg] \nonumber \\
&=& \EE^y\bigg[\int\limits_0^{\tau_D}f(Y(t))dt+g(Y({\tau_D}))\bigg]
    =J^{\tau_D}(y)\leq\Phi(y) \;.
\end{eqnarray}
Combining \eqref{4.9} and \eqref{4.11}, we get
$$
\phi(y)\geq\Phi(y)\geq J^{\tau_D}(y)=\phi(y).
$$
We conclude that
$$\phi(y)=\Phi(y)\text{ for all } \,y\in S\;.$$
Moreover, the stopping time $\widehat{\tau}$ defined by
$$
\widehat{\tau}(y,\omega)=
\begin{cases}
0 \text{ for } y\notin \D,\\
\tau_D \text{ for } y\in \D,
\end{cases}
$$
is optimal. We conclude that $\tau_{\D}$ is also optimal.
\fproof

\section{Applications}
In this section we illustrate our main result by applying it to solve explicitly some optimal stopping problems involving conditional McKean-Vlasov equations.
\subsection{The optimal time to sell}
Suppose we have a property, e.g. a house, with value $X(t)$ satisfying a McKean-Vlasov jump
diffusion with common noise given by
\small
\begin{equation}
\left\{ 
\begin{array}{ll}
dX(t)= & \mathbb{E}\left[ X(t)\mid \mathcal{F}_{t}^{(1)}\right] \Big(\alpha
_{0}(t)dt+\sigma _{1}(t)dB_{1}(t)+\sigma _{2}(t)dB_{2}(t) +\int_{\mathbb{R}}\gamma _{0}(t,\zeta )\widetilde{N}(dt,d\zeta
)\Big), \\ 
X(s)= & x>0; \quad t\geq s,
\end{array}
\right.   \label{(1)}
\end{equation}%
where $\alpha _{0},\sigma _{1},\sigma _{2}$ are positive constants and $%
\gamma _{0}$ satisfies $-1<\gamma _{0}(t,\zeta )\leq 0$. This implies that the jumps of $X$ can only be negative.\\
Suppose that if we sell the investment at a time $s+\tau $, where $\tau \in
\mathcal{T} $, the expected net profit is given by
\begin{equation}
J^{\tau }(s,x,\mu )=\mathbb{E}^{(s,x,\mu ) }\left[ e^{-\rho (s+\tau )}\Big(\EE\left[
X(\tau)\mid \mathcal{F}_{\tau }^{(1)}\right] -a\Big)\ind{} _{\tau <\infty }\right] ,
\label{(2)}
\end{equation}%
where $\rho >0$ represents the discounting coefficient and $a>0$ is the tax. We want to
find the optimal time $\tau ^{\ast }$ to sell, i.e. to find $\tau ^{\ast
}\in \mathcal{T} $ and $\Phi (s,x,\mu )$, such that%
\begin{equation}
J^{\tau ^{\ast }}(s,x,\mu )=\sup_{\tau \in \mathcal{T} }J^{\tau }(s,x,\mu
)=:\Phi (s,x,\mu ).  \label{(3)}
\end{equation}%
Comparing with Theorem 3.2, we see that in this case we have $d=1,m=2,k=1$ and%
\begin{equation}
\alpha _{1}=\alpha _{0}\left\langle \mu ,q\right\rangle ,\beta _{1}=\sigma
_{1}\left\langle \mu ,q\right\rangle ,\beta _{2}=\sigma _{2}\left\langle \mu
,q\right\rangle ,\gamma (s,x,\mu ,\zeta )=\gamma _{0}(t,\zeta )\left\langle
\mu ,q\right\rangle .  \label{(4)}
\end{equation}%
Here we have put $q(x)=x$ so that $\left\langle \mu_t ,q\right\rangle =\EE\left[ X(t)\mid 
\mathcal{F}_{t}^{(1)}\right] .$ Therefore the operator $\mathcal{G}$ takes the form:
\begin{eqnarray}
\mathcal{G}{\varphi }(s,x,\mu ) &=&\frac{\partial \varphi }{\partial s}+\alpha
_{0}\left\langle \mu ,q\right\rangle \frac{\partial \varphi }{\partial x}%
+\left\langle \nabla _{\mu }\varphi ,A_{0}^{\ast }\mu \right\rangle 
\label{(5)} \\
&&+\frac{1}{2}(\sigma _{1}^{2}+\sigma _{2}^{2})\left\langle \mu
,q\right\rangle ^{2}\frac{\partial ^{2}\varphi }{\partial x^{2}}+\frac{1}{2}%
\sigma _{1}\left\langle \mu ,q\right\rangle \frac{\partial }{\partial x}%
\left\langle \nabla _{\mu }\varphi ,A_{1}^{\ast }\mu \right\rangle  
\nonumber \\
&&+\frac{1}{2}\left\langle A_{1}^{\ast }\mu ,\left\langle D_{\mu
}^{2}\varphi ,A_{1}^{\ast }\mu \right\rangle \right\rangle   \nonumber \\
&&+\int_{\mathbb{R}}\left\{ \varphi (s,x+\gamma _{0}\left\langle \mu
,q\right\rangle ,\mu )-\varphi (s,x,\mu )-\gamma _{0}\left\langle \mu
,q\right\rangle \frac{\partial }{\partial x}\varphi (s,x,\mu )\right\} \nu
(d\zeta ),  \nonumber
\end{eqnarray}%
where%
\begin{equation}
A_{0}^{\ast }\mu =-D[\alpha _{0}\left\langle \mu ,q\right\rangle \mu ]+\frac{%
1}{2}D^{2}[(\sigma _{1}^{2}+\sigma _{2}^{2})\left\langle \mu ,q\right\rangle
^{2}\mu ],  \label{(6)}
\end{equation}%
and%
\begin{equation}
A_{1}^{\ast }\mu =-D[\sigma _{1}\left\langle \mu ,q\right\rangle \mu ].
\label{(7)}
\end{equation}%
The adjoints of the last two operators are%
\begin{equation}
A_{0}\mu =\alpha _{0}\left\langle \mu ,q\right\rangle D\mu +\frac{1}{2}%
[(\sigma _{1}^{2}+\sigma _{2}^{2})\left\langle \mu ,q\right\rangle
^{2}D^{2}\mu   \label{(8)}
\end{equation}%
and%
\begin{equation}
A_{1}\mu =\sigma _{1}\left\langle \mu ,q\right\rangle D\mu .  \label{(9)}
\end{equation}%
As a candidate for the value function $\Phi $ let us try a function of the
form%
\begin{equation}
\varphi (s,x,\mu )=\left\{ 
\begin{array}{c}
\psi(s)\left\langle \mu
,q\right\rangle ^{\lambda } ;\qquad\left\langle \mu ,q\right\rangle \in D,\\ 
e^{-\rho s}(\left\langle \mu ,q\right\rangle -a) ;\qquad\left\langle \mu
,q\right\rangle \notin D, 
\end{array}%
\right.   \label{(10)}
\end{equation}%
where $\lambda \in \mathbb{R}$ is a constant and $\psi \left( \cdot \right) $ is
a deterministic function, both to be determined. Note that, by the chain rule (see Appendix), we get
\begin{equation}
(\nabla _{\mu }\left\langle \mu ,q\right\rangle ^{\lambda })(h)=\lambda
\left\langle \mu ,q\right\rangle ^{\lambda -1}\left\langle h,q\right\rangle ;
\label{(11)}
\end{equation}%
and%
\begin{equation}
D^{2}\left\langle \mu ,q\right\rangle ^{\lambda }(h,k)=\lambda (\lambda
-1)\left\langle \mu ,q\right\rangle ^{\lambda -2}\left\langle
h,q\right\rangle \left\langle k,q\right\rangle.   \label{(12)}
\end{equation}%
Therefore
\begin{eqnarray}
\left\langle \nabla _{\mu }\varphi ,A_{0}^{\ast }\mu \right\rangle  &=&\psi
(s)\lambda \left\langle \mu ,q\right\rangle ^{\lambda -1}\left\langle
A_{0}^{\ast }\mu ,q\right\rangle  =\psi (s)\lambda \left\langle \mu ,q\right\rangle ^{\lambda
-1}\left\langle \mu ,A_{0}q \right\rangle \label{(13)} \\
&=&\psi (s)\lambda \left\langle
\mu ,q\right\rangle ^{\lambda -1}\alpha _{0}\left\langle \mu ,q\right\rangle 
=\alpha _{0}\psi (s)\lambda \left\langle \mu ,q\right\rangle ^{{\lambda }},
\nonumber
\end{eqnarray}
and%
\begin{eqnarray}
\frac{1}{2}\left\langle A_{1}^{\ast }\mu ,\left\langle D_{\mu}^{2} \varphi
,A_{1}^{\ast }\mu \right\rangle \right\rangle &=&\frac{1}{2}\psi (s)\lambda (\lambda -1)\left\langle \mu ,q\right\rangle
^{\lambda -2}\left\langle A_{1}^{\ast }\mu ,q\right\rangle \left\langle
A_{1}^{\ast }\mu ,q\right\rangle    \label{(14)} \\
&=&\frac{1}{2}\psi (s)\lambda (\lambda -1)\left\langle \mu ,q\right\rangle
^{\lambda -2}\left\langle \mu ,A_1 q\right\rangle \left\langle \mu
,A_1 q\right\rangle   \nonumber \\
&=&\frac{1}{2}\psi (s)\lambda (\lambda -1)\left\langle \mu ,q\right\rangle
^{\lambda -2}\sigma _{1}^{2}\left\langle \mu ,q\right\rangle ^{2}  \nonumber
\\
&=&\frac{1}{2}\psi (s)\lambda (\lambda -1)\sigma _{1}^{2}\left\langle \mu
,q\right\rangle ^{\lambda }.  \nonumber
\end{eqnarray}%
Moneover, note that%
\begin{equation}
\int_{\mathbb{R}}\left\{ \varphi (s,x+\gamma _{0}\left\langle \mu
,q\right\rangle ,\mu )-\varphi (s,x,\mu )-\gamma _{0}\left\langle \mu
,q\right\rangle \frac{\partial \varphi }{\partial x}(s,x,\mu )\right\} \nu
(d\zeta )=0,  \label{(15)}
\end{equation}
since $\varphi $ does not depend on $x.$ Substituting (\ref{(13)}), (\ref{(14)}) and (\ref{(15)}) into (\ref{(5)}), we get
\begin{equation}
\mathcal{G}\varphi (s,x,\mu )=\left\langle \mu ,q\right\rangle ^{\lambda }[\psi
'(s)+\psi (s)\{\alpha _{0}\lambda +\frac{1}{2}\sigma _{1}^{2}\lambda
(\lambda -1)\}].  \label{(16)}
\end{equation}
Let us guess that the continuation region has the form%
\begin{equation}
\D=\{(s,x,\mu );\left\langle \mu ,q\right\rangle <\xi \},  \label{(17)}
\end{equation}
for some $\xi >a$ (to be determined).\\
In our case we have
\begin{equation}
f=0\text{ and }g(s,x,\mu )=e^{-\rho s}(\left\langle \mu ,q\right\rangle -a).
\label{(18)}
\end{equation}
Therefore we are required to have
\begin{equation}
\mathcal{G}\varphi (s,x,\mu )=\left\langle \mu ,q\right\rangle ^{\lambda }[\psi
'(s)+\psi (s)\{\alpha _{0}\lambda +\frac{1}{2}\sigma _{1}^{2}\lambda
(\lambda -1)\}]=0\text{ for }\left\langle \mu ,q\right\rangle <\xi, 
\label{(19)}
\end{equation}
and
\begin{equation}
\mathcal{G}\varphi (s,x,\mu )=e^{-\rho s}(\left\langle \mu ,q\right\rangle -a)\text{
when }\left\langle \mu ,q\right\rangle =\xi.   \label{(20)}
\end{equation}
From (\ref{(19)}), we get
\begin{equation}
\psi (s)=\psi (0)\exp (-(\alpha _{0}\lambda +\frac{1}{2}\sigma
_{1}^{2}\lambda (\lambda -1))s), \label{(21)}
\end{equation}
for some constant $\psi (0).$\\
By \eqref{(10)} we see that continuity of $\varphi $ at $\partial \D$ implies that
\begin{equation}
\xi ^{\lambda }\psi (s)=e^{-\rho s}(\xi -a),\text{ i.e. } \psi (s)=\frac{%
\xi -a}{\xi ^{\lambda }}e^{-\rho s}.  \label{(22)}
\end{equation}
Hence
\begin{equation}
\psi (0)=\frac{\xi -a}{\xi ^{\lambda }},  \label{(23)}
\end{equation}
and, by (\ref{(21)}) and (\ref{(22)}),
\begin{equation}
\psi (0)\exp (-(\alpha _{0}\lambda +\frac{1}{2}\sigma _{1}^{2}\lambda
(\lambda -1))s)=\psi (s)=\psi (0)e^{-\rho s}.  \label{(24)}
\end{equation}
Therefore,
\begin{equation}
\alpha _{0}\lambda +\frac{1}{2}\sigma _{1}^{2}\lambda (\lambda -1)=\rho.
\label{(25)}
\end{equation}
Equation (\ref{(25)}) gives two possible values of $\lambda :$%
\begin{equation}
\lambda _{i}=\sigma _{1}^{-2}\left[ \frac{1}{2}\sigma _{1}^{2}-\alpha_0 \pm 
\sqrt{(\alpha_0 -\frac{1}{2}\sigma _{1}^{2})^{2}+2\rho \sigma _{1}^{2}}\right]
;\text{ }i=1,2,  \label{(26)}
\end{equation}
where $\lambda _{2}<0<\lambda _{1}.$\\
Since the value function must be bounded for
$\left\langle \mu ,q\right\rangle \leq \xi ,$ we choose $\lambda =\lambda_{1}$ and get
\begin{equation}
\varphi (s,x,\mu )=\psi (s)\left\langle \mu ,q\right\rangle ^{\lambda _{1}}=%
\frac{\xi -a}{\xi ^{\lambda }}e^{-\rho s}\left\langle \mu
,q\right\rangle ^{\lambda _{1}};\left\langle \mu ,q\right\rangle \leq \xi.
\label{(27)}
\end{equation}%
If we assume that $\alpha <\rho ,$ then we see after some computation that $\gamma _{1}>1$. \\

It remains to find the value of $\xi $ which makes the function $\varphi $
defined on all $\mathbb{R\times R\times M}$ by 
\begin{equation}
\varphi (s,x,\mu )=\left\{ 
\begin{array}{c}
\frac{\xi -a}{\xi ^{\lambda_1 }}e^{-\rho s}\left\langle \mu
,q\right\rangle ^{\lambda _{1}} ;\qquad\left\langle \mu ,q\right\rangle \leq\xi,  \\ 
e^{-\rho s}(\left\langle \mu ,q\right\rangle -a);\qquad \left\langle \mu
,q\right\rangle \geq \xi, 
\end{array}
\right.   \label{(30)}
\end{equation}
continuously differentiable at $\left\langle \mu ,q\right\rangle =\xi.$\\
Note that
\[
\lim_{\left\langle \mu ,q\right\rangle \rightarrow \xi ^{-}}\nabla _{\mu
}\varphi(s,x,\mu)= \frac{\xi - a}{\xi ^{\lambda_1}} \lambda_1 \xi^{\lambda_1} e^{-{\rho s}}=\lambda_1 (\xi - a)e^{-\rho s}, 
\]%
and 
\[
\lim_{\left\langle \mu ,q\right\rangle \rightarrow \xi ^{+}}\nabla _{\mu
}(e^{-\rho s}(\left\langle \mu ,q\right\rangle -a))=e^{-\rho s}\xi . 
\]%
Hence for differentiability at $\left\langle \mu ,q\right\rangle =\xi $ we
should have
$\lambda _{1}(\xi -a)=\xi ,$ or
\begin{equation}
\xi =\xi ^{\ast }:=\frac{\lambda _{1}a}{\lambda _{1}-1},  \label{(31)}
\end{equation}%
with this choice of $\xi =\xi ^{\ast }$ it is easy to verify that the
function $\varphi $ defined by \eqref{(30)}, \eqref{(31)} satisfies all the requirements
of the verification theorem, and we have proved the following:

\begin{theorem}
Assume that $\alpha < \rho$. Then the  value function $\Phi $ of the optimal stopping problem (\ref{(3)}) is
given by $\Phi =\varphi ,$
where $\varphi $ is defined in (\ref{(30)}), with $\xi =\xi ^{\ast }$
given by (\ref{(31)}).\\

An optimal stopping time is $\tau^{*}:=\tau_{\D}$, where $\D$ is given by \eqref{(17)} 
with $\xi=\xi^{*}$, i.e.
\begin{align}
\tau^{*}=\inf\{t>0; \EE\left[ X(t)\mid 
\mathcal{F}_{t}^{(1)}\right] \geq \xi^{*} \}.
\end{align}
\end{theorem}

\subsection{Optimal time to quit a project}
Suppose the value $X(t)$ at time $t$ of a resource is given by the equation 
\begin{equation}
dX(t)=\sigma _{1}dB_{1}(t)+\sigma _{2}dB_{2}(t)+\int_{\mathbb{R}}\gamma
_{0}(\zeta )\widetilde{N}(dt,d\zeta ), \quad X(s)=x,  \label{1a}
\end{equation}%
and suppose that if we stay in the project till time $\tau$ the expected discounted profit is%
\begin{equation}
J^{\tau }(s,x,\mu )=\mathbb{E}^{(s,x,\mu) }\left[ \int_{0}^{\tau }e^{-\rho
(s+t)}\mathbb{E}\left[ X(t)\mid \mathcal{F}_{t}^{(1)}\right] dt\right],
\label{2a}
\end{equation}
where $\rho > 0$ is a given discounting coefficient, $\sigma _{1},\sigma _{2}$ and $\gamma_0$ are given nonzero constants.

The problem is to find the stopping time $\tau^{*}$ which maximizes this expected net profit.
This is an optimal stopping problem of the type discussed in Section 3. Here we have 

$d=1,m=2,\alpha =0,\beta _{11}=\sigma _{1},\beta _{12}=\sigma _{2},\gamma
(s,x,\mu ,\zeta )=\gamma _{0}(\zeta ), f(s,x,\mu )=e^{-\rho s}\left\langle \mu ,q\right\rangle $ and $g=0.$

\noindent If we as before put $q(x)=x$, then $\left\langle \mu _{t},q\right\rangle =\mathbb{E}\left[
X(t)\mid \mathcal{F}_{t}^{(1)}\right]$. 
In this case we get
\small
\begin{eqnarray}
\mathcal{G}\varphi (s,x,\mu ) &=&\frac{\partial \varphi }{\partial s}+\left\langle
\nabla _{\mu }\varphi ,A_{0}^{\ast }\mu \right\rangle   +\frac{1}{2}(\sigma _{1}^{2}+\sigma _{2}^{2})\frac{\partial ^{2}\varphi }{%
\partial x^{2}}+\frac{1}{2}\sigma _{1}\frac{\partial }{\partial x}%
\left\langle \nabla _{\mu }\varphi ,A_{1}^{\ast }\mu \right\rangle  
\label{3a}\\
&&+\frac{1}{2}\left\langle A_{1}^{\ast }\mu ,\left\langle D_{\mu
}^{2}\varphi ,A_{1}^{\ast }\mu \right\rangle \right\rangle   \nonumber \\
&&+\int_{\mathbb{R}}\left\{ \varphi (s,x+\gamma _{0},\mu )-\varphi (s,x,\mu
)-\gamma _{0}\frac{\partial }{\partial x}\varphi (s,x,\mu )\right\} \nu
(d\zeta ),  \nonumber
\end{eqnarray}%
where%
\begin{equation}
A_{0}^{\ast }\mu =\frac{1}{2}D^{2}[(\sigma _{1}^{2}+\sigma _{2}^{2})\mu
],\text{ which is the adjoint of } A_{0}\mu =\frac{1}{2}(\sigma _{1}^{2}+\sigma _{2}^{2})D^{2}\mu, 
\label{4a}
\end{equation}%
and%
\begin{equation}
A_{1}^{\ast }\mu =-D[\sigma _{1}\mu ],\text{ which is the adjoint of } A_{1}\mu =\sigma _{1}D\mu .
\label{5a}
\end{equation}%
We guess that the continuation region has the form%
\begin{equation}
\D=\{(s,x,\mu );\left\langle \mu ,q\right\rangle > \eta \},  \label{6a}
\end{equation}%
for some $\eta <0$ ($t_{0}$ be determined) and that the value function is of
the form%
\begin{equation}
\varphi (s,x,\mu )=\left\{ 
\begin{array}{l}
e^{-\rho s}F(\left\langle \mu ,q\right\rangle )\text{\qquad for }%
\left\langle \mu ,q\right\rangle \geq \eta,  \\ 
0\text{\qquad \qquad \qquad \qquad for}\left\langle \mu ,q\right\rangle <\eta, 
\end{array}%
\right.   \label{7a}
\end{equation}%
for some $C^{1}$ function $F$ ($t_{0}$ be determined).
Then by the chain rule,%
\[
\nabla _{\mu }\varphi (h)=e^{-\rho s}F^{\prime }(\left\langle \mu
,q\right\rangle )\left\langle h,q\right\rangle,
\]%
so that%
\begin{eqnarray}
\left\langle \nabla _{\mu }\varphi ,A_{0}^{\ast }\mu \right\rangle 
&=&e^{-\rho s}F^{\prime }(\left\langle \mu ,q\right\rangle )\left\langle
A_{0}^{\ast }\mu ,q\right\rangle   \label{8a} \\
&=&e^{-\rho s}F^{\prime }(\left\langle \mu ,q\right\rangle )\left\langle \mu
,A_{0}q\right\rangle =0.  \nonumber
\end{eqnarray}%
Similarly, by the chain rule,%
\[
D_{\mu }^{2}\varphi (h,k)=F^{\prime \prime }(\left\langle \mu
,q\right\rangle )\left\langle h,q\right\rangle \left\langle k,q\right\rangle
e^{-\rho s},
\]%
and hence 
\begin{eqnarray}
\left\langle A^{\ast },\mu ,\left\langle D_{\mu }^{2}\varphi ,A^{\ast }\mu
\right\rangle \right\rangle  &=&F^{\prime \prime }(\left\langle \mu
,q\right\rangle )\left\langle A^{\ast },\mu ,q\right\rangle \left\langle
A^{\ast },\mu ,q\right\rangle e^{-\rho s}  \label{9a} \\
&=&F^{\prime \prime }(\left\langle \mu ,q\right\rangle )\left\langle \mu
,A,q\right\rangle \left\langle \mu ,A,q\right\rangle e^{-\rho s}=F^{\prime
\prime }(\left\langle \mu ,q\right\rangle )\sigma _{1}^{2}e^{-\rho s} .
\nonumber
\end{eqnarray}%
Substituting this into \eqref{3a} we get that, for $(s,x,\mu )\in \D,$%
\begin{equation}
\mathcal{G}\varphi (s,x,\mu )+f(s,x,\mu )=e^{-\rho s}\left[ -\rho F(\left\langle
\mu ,q\right\rangle )+\frac{1}{2}\sigma _{1}^{2}F^{^{\prime \prime
}}(\left\langle \mu ,q\right\rangle )+\left\langle \mu ,q\right\rangle %
\right] =0.  \label{10a}
\end{equation}%
Hence the function $F$ must satisfy%
\begin{equation}
z-\rho F(z)+\frac{1}{2}\sigma _{1}^{2}F^{^{\prime \prime }}(z)=0,\text{
where }z=\left\langle \mu ,q\right\rangle .  \label{11a}
\end{equation}%
This equation has the general solution%
\begin{equation}
F(z)=\frac{1}{\rho }z+C_{1}e^{-\lambda z}+C_{2}e^{\lambda z}\text{ in }D,
\label{12a}
\end{equation}%
where $C_{1},C_{2}$ are constants and%
\begin{equation}
\lambda =\sqrt{\frac{2\rho }{\sigma _{1}^{2}}}.  \label{13a}
\end{equation}%
Since we do not expect that $\Phi $ will grow exponentially when $%
z\rightarrow \infty ,$ we guess that $C_{2}=0,$ and get by \eqref{7a}%
\begin{equation}
\varphi (s,x,\mu )=\left\{ 
\begin{array}{l}
e^{-\rho s}\left[ \frac{1}{\rho }\left\langle \mu ,q\right\rangle
+C_{1}e^{-\lambda \left\langle \mu ,q\right\rangle }\right] ;\left\langle
\mu ,q\right\rangle \geq \eta,  \\ 
0\text{ for }\left\langle \mu ,q\right\rangle <\eta. 
\end{array}%
\right.   \label{14a}
\end{equation}%
If we require continuity of $\varphi $ at $\left\langle \mu ,q\right\rangle =\eta $ we get 
the equation%
\begin{equation}
\frac{1}{\rho }\eta +C_{1}e^{-\lambda \eta }=0.  \label{15a}
\end{equation}%
Requiring that $\varphi $ is $e^{1}$ at $\left\langle \mu ,q\right\rangle
=\eta $ gives the equation
\begin{equation}
\frac{1}{\rho }-\lambda C_{1}e^{-\lambda \eta }=0.  \label{16a}
\end{equation}%
Combining \eqref{15a} and \eqref{16a} we get%
\begin{equation}
\eta =\eta ^{\ast }=-\frac{1}{\lambda }=-\sqrt{\frac{2\rho }{\sigma _{1}^{2}}%
}  \label{17a}
\end{equation}%
and%
\begin{equation}
C_{1}=-\frac{\eta ^{\ast }}{\rho }e^{\lambda \eta ^{\ast }}.  \label{(18)}
\end{equation}%
We can now easily verify that with this choice of $\eta =\eta ^{\ast }$ and $%
C_{1}$ given by \eqref{17a} and \eqref{(18)} the function $\varphi $
defined by \eqref{14a} satisfies all the condition of are verification
theorem.\\
Therefore we have proved the following theorem:

\begin{theorem}
The value function $\Phi $ of the optimal stopping problem%
\[
\Phi (s,x,\mu )=\mathbb{E}^{(s,x,\mu) }\left[ \int_{0}^{\tau }e^{-\rho (s+t)}%
\mathbb{E}\left[ X(t)\mid \mathcal{F}_{t}^{(1)}\right] dt\right] 
\]%
is given by
\[
\Phi (s,x,\mu )=\left\{ 
\begin{array}{l}
e^{-\rho s}\left[ \frac{1}{\rho }\left\langle \mu ,q\right\rangle -\frac{%
\eta ^{\ast }}{\rho }e^{\eta ^{\ast }(\left\langle \mu ,q\right\rangle -\eta
^{\ast })}\right] ;\text{for }\left\langle \mu ,q\right\rangle \geq \eta
^{\ast } \\ 
0\text{ for }\left\langle \mu ,q\right\rangle <\eta ^{\ast },
\end{array}
\right. 
\]%
where%
\[
\eta ^{\ast }=-\sqrt{\frac{z\rho }{\sigma _{1}^{2}}.}
\]%
The corresponding optimal stopping time is the first time the conditional density of $X(t)$  hits or goes below the threshold $\eta^{*}$, i.e.
\[
\tau ^{\ast }=\inf \Big\{t>0;\mathbb{E}\left[ X(t)\mid \mathcal{F}_{t}^{(1)}%
\right] \leq \eta ^{\ast }\Big\}.
\]
\end{theorem}

\section{Appendix: Double Fr\' echet derivatives}
In this section we recall some basic facts we are using about the Fr\' echet derivatives of a function $f: V \mapsto W$, where $V,W$ are given Banach spaces.
\begin{definition}
We say that $f$ has a Fr\' echet derivative $\nabla_xf=Df(x)$ at $x \in V$ if there exists a bounded linear map $A:V \mapsto W$ such that
\begin{align*}
\lim_{h \rightarrow 0} \frac{||f(x+h)-f(x)-A(h)||_{W}}{||h||_{V}} =0.
\end{align*}
Then we call $A$ the Fr\' echet derivative of $f$ at x and we put $Df(x) =A$.
\end{definition}
Note that $Df(x) \in L(V,W)$ (the space of bounded linear functions from $V$ to $W$),  for each $x$.

\begin{definition}
We say that $f$ has a double Fr\' echet derivative $D^2 f(x)$ at $x$ if there exists a bounded bilinear map $A(h,k): V \times V \mapsto W$ such that
\begin{align*}
\lim_{k \rightarrow 0} \frac{||Df(x+k)(h)-Df(x)(h)-A(h,k)||_{W}}{||h||_{V}} =0.
\end{align*}
\end{definition}

\begin{example}
\begin{itemize}
\item
Suppose $f:\mathbb{M} \mapsto \mathbb{R}$ is given by
\begin{align*}
f(\mu)=\<\mu,q\>^2, \text{ where } q(x)=x.
\end{align*}
Then 
\begin{align*}
f(\mu +h) -f(\mu)&= \<\mu+h,q\>^2 -\<\mu,q\>^2\nonumber\\
&= 2 \<\mu,q\> \<h,q\>+\<h,q\>^2, 
\end{align*}
so we see that
\begin{equation*}
Df(\mu)(h)=2\<\mu,q\>\<h,q\>.
\end{equation*}
To find the double derivative we consider
\begin{align*}
&Df(\mu+k)(h)-Df(\mu)(h)\nonumber\\
&=2\<\mu+k,q\>\<h,q\>-2\<\mu,q\>\<h,q\>\nonumber\\
&=2\<k,q\>\<h,q\>,
\end{align*}
and we conclude that
\begin{equation*}
D^2f(\mu)(h,k)=2\<k,q\>\<h,q\>.
\end{equation*}
\item
Next assume that $g:\mathbb{M} \mapsto \RR$ is given by $g(\mu)=\<\mu,q\>$.
Then, proceeding as above we find that
\begin{align*}
Dg(\mu)(h)&=\<h,q\> \text{ (independent of } \mu)\\
&\text{ and }\nonumber\\
D^2g(\mu) &=0.
\end{align*}

\end{itemize}
\end{example}


\end{document}